# Mixed Mode Oscillations and Bifurcation Mechanism in a Nonlinear Beam-Elastic Foundation Under Parametric and External Excitations

Sobhan Mohammadi, Keegan J. Moore[*]

Daniel Guggenheim School of Aerospace Engineering, Georgia Institute of Technology, Atlanta, GA 30332, United States

*Corresponding author: (K.J. Moore)

E-mail address: kmoore@gatech.edu

## ABSTRACT

This paper aims to study existence condition of possible bursting oscillations generated by low frequency excitation of a nonlinear vibratory system in the presence of parametric excitation. Slow-fast dissection technique and numerical bifurcation analysis are employed to extract qualitative changes in system response originated from its nonlinear dynamics. Role of all parameters of elastic foundation and excitation model are studies and it is shown that the system exhibits the phenomena of folding, cusp and Bogdanov-Takens bifurcations which are potentially routes to bi-stability and chaos. It can be found that slow excitation of the nonlinear foundation is the main generating factor of fold bifurcation and stiffness of elastic foundation has a remarkable effect on stability region of the beam. In addition, the base excitation of an elastic foundation in form of a traveling wave, adds multi-frequency excitation and parametric resonances necessarily to the system. This study showed investigating nonlinear oscillator under low frequency excitations in framework of slow-fast plays an invaluable role in understanding instabilities in systems that are not captured by standard methods.

**Keywords:** Mixed-mode oscillations, slow-fast dynamic, nonlinear oscillations, parametric excitations

## 1. Introduction

Dynamical systems with multiple time scales are abundant in science and engineering. Among the complex behaviors observed in such systems, mixed-mode oscillations (MMOs) also referred to as bursting oscillations—are particularly intriguing. MMOs typically arise in slow-fast dynamical systems, where the interaction between slow and fast subsystems leads to the emergence of oscillatory patterns that combine small-amplitude fluctuations with intermittent large spikes. These patterns reflect alternating phases of gradual evolution and rapid transitions in the system's state [1]. Bursting oscillations were classified and introduced for first time by Rinzel [2, 3] in field of neural dynamics, where bursting oscillations were found to be caused by the transition of system between different resting and active states. In Other words, in bursting we expect to see a series of spiking amplitudes separated by a period of inactivity. Rinzel also introduced the slow-fast technique (also known as slow-fast dissection) [3] to analyze and explain this strongly nonlinear behavior. Building on this research, the study of systems with bursting oscillations has attracted researchers from many different areas. Zhou and coworkers [4] studied bursting oscillations of a chemical reacting system under two slow parametric excitations and found four different bursting patterns that emerge in that system. Zhang and Bi [5] investigated the occurrence of bursting oscillations in a dynamical system in presence of Hopf bifurcation using the method of normal forms. They revealed the existence of a double Hopf bifurcation and found that increasing of amplitude of the excitation led to switching of the trajectory from a single-mode to mixed-mode branch. Ma et al. [6] applied the slow-fast method to investigate four compound bursting behaviors in a Mathieu-Van der Pol-Duffing oscillator. They studied possible routes to compound bursting dynamics and found that these bursting behaviors originated through interactions of attractors and bifurcation structures. Ma et al. [7] reported four complex bursting patterns in cubic-quintic Duffing-Van der Pol system under two external excitations and argued that amplitude of excitation has a key role in system response.

Bursting oscillations have also been observed in non-smooth dynamical systems, systems that are discontinuous in nature such as freeplay and impact in mechanical devices or switching elements in electronic devices. Non-smooth dynamical systems could be classified in terms of degree of smoothness. Filipov is an important subclass of non-smooth systems with degree of smoothness equal to one, which means the flows in system are separated by one piecewise smooth boundary. According to the nature of the boundary and implied discontinuity, there is possibility of occurring four types of bifurcations specific to this class of systems known as crossing-sliding, grazing-sliding, switching-sliding and adding-sliding. Interested readers are referred to the Ref.[8]. Han and Bi [9] studied the sliding dynamics of a Duffing oscillator under slow excitation with a switching property. They found six patterns of slow-fast oscillations that can occur and that the switching threshold frequency has a significant effect on system behavior. They also demonstrated that the phenomena could be reproduced using a similar piecewise-smooth dynamical system. In a similar work, Xu et al. [10] investigated sliding bursting for a non-smooth coupled Duffing-Van der Pol oscillator under pure single external excitation. *Sliding* in non-smooth dynamic systems context means one of the system variables remains constant as the system passes through a switching manifold. Subsequently, *sliding bursting* implies bursting behavior which occurs during sliding. They reported that Hopf-like bifurcations facilitate the transition between spiking and quiescent states. Wang and coworkers [11] studied bursting oscillations of a piecewise-smooth electrical circuit in the presence of a codimension-2 bifurcation. They concluded that the transition between spiking and quiescent states is governed by limit cycle bifurcation. Zhao et al. [12] also addressed noticeable slow-fast dynamics for a slowly excited Duffing oscillator with ascending-descending switching excitations. They found two bifurcation structures and several novel slow-fast dynamics including jumps hysteresis loop and switching frequency failure phenomenon.

Li and Hou [13] studied generating mechanism of bursting phenomenon in a piecewise linear mechanical oscillator with two time scales under stiffness perturbation and slow excitation. Using slow-fast analysis, they found that the focus-type periodic bursting oscillation could evolve in the stable response of the system. This type of bursting oscillation occurs intermittently in piecewise systems during the quiescent phase and exhibits behavior resembling a focus-type fixed point. Depending on whether the trajectory spirals inward or outward, the oscillation can be classified as stable or unstable. Furthermore, they concluded that although damping ratio and excitation amplitude regulate the system response, the bursting oscillations have unique generation mechanisms. Recently, Wei and Han [14] investigated the dynamics of sharp transitions caused by instantaneous changes due to a square wave excitation in the Rayleigh oscillator. They found that catastrophic jumps can be expected in response because of abrupt transitions caused by step excitations. They believed that observed jumps are independent of specific system and they are induced by square wave excitation. Moreover, Qian et al. [15] studied bursting dynamics in a hybrid Rayleigh-Van der Pol-Duffing oscillator under two excitations. They reported that the main bursting pattern caused by a single excitation completely changed with the addition of a second excitation. Furthermore, Zhao and coworkers [16] captured chaotic and periodic bursting patterns caused by transient chaos in a smooth three-dimensional dynamic systems and argued that in addition to the bifurcation structure, delay can also affect the bursting pattern.

Amongst the research conducted on bursting oscillations in mechanical systems, Liu and coworkers [17] applied the slow-fast method to a pendulum with an irrational nonlinearity with parametric and external excitation. They considered one, two and three parameter bifurcations and found diverse bursting dynamics including hysteresis cycle and jumping. Among works done in structural stability area with focus on slow-fast dynamics, Hao et al.[18] explored bursting oscillations in the coupled bending-torsions of a sandwich conical panel using first modes for bending and torsion and captured pitchfork and symmetry breaking bifurcations. More interestingly, they found that bursting oscillations only occur around static equilibriums at very small amplitudes. They concluded that preload has a meaningful effect on dynamic response of the panel as well as the torsional response undergo unpredictable large amplitudes.

Slow-fast dynamical system encompasses many broad classes of problems in different areas, hence bursting oscillations are ubiquitous. Some research conducted on slow-fast systems could be found in electrical circuits (Lin et al. [19], Han et al. [20], Bao et al. [21]), ecology (Li et al. [22]), hydraulic power turbines (Zhang et al. [23] and Li et al. [24]), and neural dynamics (Wang et al. [25], Hoppensteadt and Izhikevich [26], Izhikevich [1], and Fallah [27]).

Research works done in slow-fast dynamics of structural systems are quite rare and this research focuses on the bursting oscillations that arise in the dynamics of a Euler-Bernoulli beam with a nonlinear, elastic foundation subjected to base excitation in the form of a traveling wave which exert forced displacement. We explore bursting patterns of a smooth slow-fast dynamical system. Stability of equilibrium points are discussed and the bifurcation analysis using MATCONT [28], as a continuation package for bifurcation studies is performed. To this aim, the forcing frequency term is taken as slow variable. Excitation frequency and amplitude of excitation are considered as bifurcation parameters and based on their variations; different bursting oscillation pattern are investigated. Two parameter bifurcation analysis were employed to explore cascades of possible transitions after folding behavior.

## 2. Mathematical Model
### 2.1. Beam Model
The dynamics of a beam supported by an elastic foundation has remarkable importance for design and modeling purposes in diverse areas covers from rail engineering to MEMS. Therefore, the interaction of beams with foundations brings new rich dynamics that need to be investigated. The equation of interest of this work is a simplified derivation from research performed with aim of modeling ship vibrations under sea wave excitation. The schematic of a beam resting on an elastic foundation with cubic hardening type springs is shown in Figure 1. The system is excited through the elastic foundation by a traveling wave $\eta(x,t)$ that moves with phase velocity $V_p$. The traveling wave exerts a force on the beam through the nonlinear foundation, such that the system is undergoes base excitation. Furthermore, due to the nonlinearity in the foundation, the system also experiences parametric resonance as will be seen in derivation of governing equation.

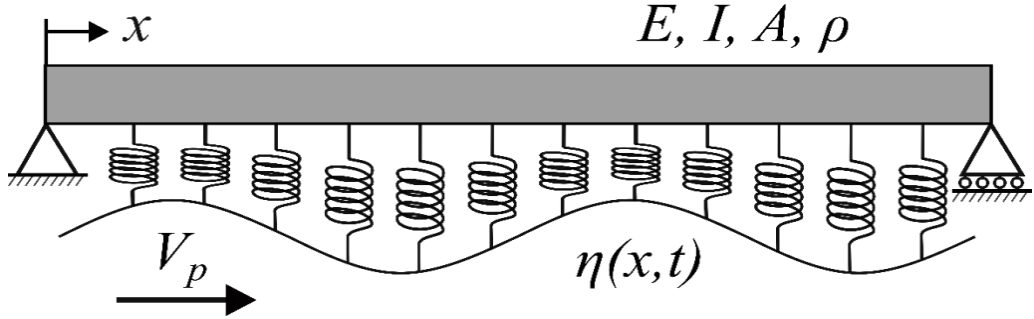

Figure 1. The beam on elastic foundation under traveling wave excitation.

The equation of motion governing the beam is

$$\rho A u_{tt} + EI u_{xxxx} + k(u - \eta(x,t)) + \alpha(u - \eta(x,t))^3 = 0, \tag{1}$$

In which $k$ and $\alpha$ are linear and nonlinear stiffness coefficients respectively. By expanding the terms

$$\rho A u_{tt} + EI u_{xxxx} + (k + 3\alpha\eta^2(x,t))u - 3\alpha\eta(x,t)u^2 + \alpha u^3 = k\eta + \alpha\eta^3, \tag{2}$$

With the aim of simplicity, the analysis will be conducted on non-dimensional equations. The non-dimensional groups are defined as follows:

$$w = \frac{u}{L}, \quad \chi = \frac{x}{L}, \quad z = \frac{\eta}{L}, \quad \tau = \frac{t}{b}, \quad b^2 = \frac{\rho A L^4}{EI}, \tag{3}$$

By using definitions in Eq. (3) and substituting in Eq. (2), the new equation obtained as,

$$w_{tt} + w_{xxxx} + (\sigma + 3\gamma z^2)w - 3\gamma z w^2 + \gamma w^3 = \sigma z + \gamma z^3, \tag{4a}$$

where

$$\sigma = \frac{kL^4}{EI}, \qquad \gamma = \frac{aL^6}{EI}. \tag{4b}$$

To extract the equation governing temporal part, Galerkin's method is used to approximate the equation according as

$$w(\chi, \tau) = \sum_{n=1}^{\infty} W(\chi) q(\tau), \tag{5}$$

where $W(\chi)$ and $q(\tau)$ are the mode shape function and modal coordinates, respectively and $n$ represents the mode number. By considering simply-supported boundary conditions shown in Figure 1 and assuming the first mode shape for approximation $w(X, \tau) = \sin(\pi X) q(t)$, the governing ODE of temporal part of the solution in Eq. (5) can be derived after calculating the integral as in Eq. (6),

$$\int_0^1 W(x)[w_{tt} + w_{xxxx} + (\sigma + 3\gamma z^2)w - 3\gamma z w^2 + \gamma w^3 - \sigma z - \gamma z^3] dx, \tag{6}$$

The final result of operation in Eq. (6) leads to the following time-variable coefficients ODE:

$$\ddot{q}(t) + \left(\pi^4 + \sigma + \frac{3\gamma \eta_0^2}{2\kappa L^2 (\pi^2 - \kappa^2)}(\kappa(\pi^2 - \kappa^2) + \pi^2 \sin(\kappa)\cos(\kappa - 2\omega t))\right) q(t) -$$

$$\left(\frac{36\pi^3 \gamma \eta_0}{L(9\pi^4 + 10\pi^2 \kappa^2 + \kappa^4)}(\cos(\omega t) + \cos(\kappa - \omega t))\right) q^2(t) + \frac{3}{4}\gamma q^3(t) =$$

$$-\frac{2\pi\sigma\eta_0}{L(\pi^2 - \kappa^2)}[\cos(\omega t) + \cos(\kappa - \omega t)] +$$

$$\frac{\pi\gamma\eta_0^3}{2L^3(9\kappa^4 - 10\pi^2\kappa^2 + \pi^4)}\begin{bmatrix}3(\pi^2 - 9\kappa^2)[\cos(\omega t) + \cos(\kappa - \omega t)] + \\ (\pi^2 - \kappa^2)[\cos(3(\kappa - \omega t)) + \cos(3\omega t)]\end{bmatrix}, \tag{7}$$

### 2.2. The Equation of Interest

Building on the derivation in the previous section, the dimensionless model of forced mechanical oscillator studied in this paper after adding a damping term is written as

$$\ddot{q}(t) + \xi \dot{q}(t) + J(t) q(t) - G(t) q^2(t) + \frac{3}{4}\gamma q^3(t) = F(t), \tag{8}$$

where $\xi$ is the damping ratio, $\gamma$ is hardening term coefficient, and $J(t)$ and $G(t)$ are time-dependent coefficients of the equation which implies parametric resonance of the system, and $q(t)$ is the displacement of the oscillator. Introducing state variables, $z_1 = q(t), z_2 = \dot{q}(t)$, Eq. (8) can be written as

$$\begin{cases} \dot{z}_1 = z_2, \\ \dot{z}_2 = -\xi z_2 - J(t)z_1 + G(t)z_1^2 - \frac{3}{4}\gamma z_1^3 + F(t). \end{cases} \quad (9)$$

Implementing slow-fast analysis necessitates transforming system in Eq. (9) into autonomous form which is discussed in next section.

### 2.3. Slow-Fast Analysis

Here, we introduce briefly the process of analyzing a system using the slow-fast dissection technique introduced by Rinzel and Lee [29] which is established based on separation of slow and fast subsystems and studying them separately. The essence of a slow-fast dynamical system is co-existence and interaction of variables with different time scales. The general form of a slow-fast dynamic system is as

$$\begin{cases} \dot{z}_1 = f(z_1, z_2), \\ \dot{z}_2 = \epsilon g(z_1, z_2), \end{cases} \quad (10)$$

in which $f(z_1, z_2)$ and $g(z_1, z_2)$ are smooth functions, and $\dot{z}_1$ and $\dot{z}_2$ are defined as fast and slow subsystems respectively. The slow-fast technique is based on transforming systems with two variables into slow-fast forms with only a single slow variable. To achieve this, we need to choose one variable to serve as the slow variable, then relate the other terms to it. In this research, $\omega$ is chosen to be significantly smaller than first natural frequency of the corresponding linear systems. Therefore, $\cos \omega t$ in forcing function is chosen to be the slow variable. In addition to the slow variable, it is also needed to find corresponding expressions for all the terms dependent on $\omega t$. To this end, we leverage the de Moivre formula, which allows one to express any multipliers of the multi-frequency in terms of the corresponding slow variable. For any real number $a$ and integer $n$, it holds that

$$(\cos a + i \sin a)^n = \cos na + i \sin na, \quad (11)$$

where $i$ is imaginary unit. By expanding the left-hand side of the Eq. (11) and equating the real parts of both sides of the equation, we can obtain an equivalent expression for $\cos(nx)$ for a specified integer value of $n$. Let $d = cos(\omega t)$, then expression of $\cos(\kappa - \omega t)$ is found to be

$$\cos(\kappa - \omega t) = \cos(\kappa)d + \sin(\kappa)\sqrt{1 - d^2}, \quad (11)$$

where the Pythagorean identity is used to deal with sine terms. The equivalent expression for $\cos(3\omega t)$ and $\cos(3(\kappa - \omega t))$ are

$$\cos(3\omega t) = 4d^3 - 3d. \tag{12}$$

$$\cos(3(\kappa - \omega t)) = \cos(3\kappa)(4d^3 - 3d) + \sin(3\kappa)(4d^2 - 1)\sqrt{1 - d^2}, \tag{13}$$

With the transformation of $d$ as the slow variable, we find the final equation as

$$\begin{cases} \dot{z}_1 = z_2 \\ \dot{z}_2 = F(d) - J(d)z_1 + G(d)z_1^2 - \frac{3}{4}\gamma z_1^3 - \xi z_2 \end{cases}, \tag{14}$$

such that

$$F(d) = -\frac{2\pi\sigma\eta_0}{L(\pi^2 - \kappa^2)}\left[d + d\cos(\kappa) + \sqrt{1 - d^2}\sin(\kappa)\right] + \frac{\pi\gamma\eta_0^3}{2L^3(9\kappa^4 - 10\pi^2\kappa^2 + \pi^4)} *$$

$$\begin{bmatrix} 3(\pi^2 - 9\kappa^2)\left[d + d\cos(\kappa) + \sqrt{1 - d^2}\sin(\kappa)\right] + \\ 2(\pi^2 - \kappa^2)\cos\left(\frac{3\kappa}{2}\right)\left[(4d^3 - 3d)\cos(3\kappa) + (4d^2 - 1)\sqrt{1 - d^2}\sin(3\kappa)\right] \end{bmatrix}, \tag{15}$$

$$J(d) = \pi^4 + \sigma + \frac{3\gamma\eta_0^2}{2\kappa L^2(\pi^2 - \kappa^2)}\left(\begin{array}{c}\kappa(\pi^2 - \kappa^2) + \\ \pi^2\sin(\kappa)\left(\cos(\kappa)(2d^2 - 1) + 2d\sqrt{1 - d^2}\sin(\kappa)\right)\end{array}\right), \tag{16}$$

$$G(d) = \frac{36\pi^3\gamma\eta_0}{L(9\pi^4 + 10\pi^2\kappa^2 + \kappa^4)}\left(d + d\cos(\kappa) + \sqrt{1 - d^2}\sin(\kappa)\right). \tag{17}$$

As a result of variable change of $d = \cos\omega t$, Eq. (15) gives forcing function as a single variable. In addition, using the expressions developed in Eq. (11) through Eq. (17), the system in Eq. (9) becomes autonomous without explicit dependence on time.

## 3. Numerical Results and Bifurcations in the Autonomous System

To demonstrate the slow-fast response and bursting oscillations that arise in the dynamics of this system, we set $\xi = 0.02, \eta_0 = 0.1L, \kappa = 1, \sigma = 120, \gamma = 150,$ and $\omega = 0.01$. These values for $\sigma$ and $\gamma$ are chosen to illustrate the diversity and richness of the dynamics governing this system. We note that stiffness of the system in Eq.(8) changes with time which prevent to have an exact constant estimation of natural frequency of the sys, but as an approximation $\omega_n \sim 10$, therefore $\omega \ll \omega_n$. The system is numerically integrated for $t \in [0,2000]$ with a step size of $5e\text{-}3$ with zero initial conditions using ODE45 in MATLAB® with default tolerances. The displacement response of the system and the phase portrait for the first and second cases are presented in Figures 2a and b, respectively. As has been shown in the, bursting oscillations emerge repeatedly while for higher excitation frequencies such a pattern of oscillation never appears, and this formed the motivation of this research to investigate the generating mechanism of these oscillations. In both cases, the system undergoes slow excitations at least one order of magnitude slower than its natural frequency. Also, the system has two center of response formed around $z_1 \sim \pm 0.5$. While the system is following the slow motion imposed by the foundation, known as quiescent phase, the beam exhibits large amplitudes of oscillations which are at higher frequencies compared to the fundamental frequency which is defined as spiking phase. Also, the system moves between the states by fast transitions which means quiescent phase in this configuration of beam-elastic foundation is a short period.

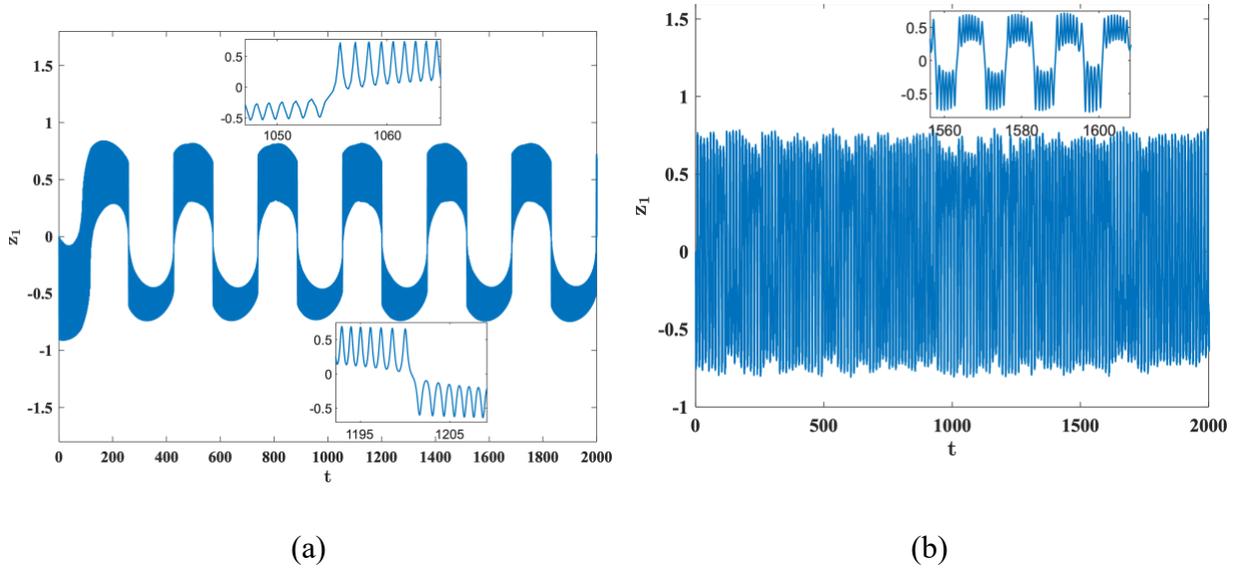

(a)          (b)

Figure 2. Bursting oscillations for $\eta_0 = 0.1, \gamma = 150, \sigma = 120$ (a) $\omega = 0.01$ (b) $\omega = 1$. Fast transition between equilibrium states could be seen in magnified subplots.

Using the slow-fast transformed system in Eq. (14), the equilibrium points of system are obtained by setting $\dot{z}_1 = 0$ and $\dot{z}_2 = 0$, leading to

$$\begin{cases} z_2 = 0, \\ -J(d)z_1 + G(d)z_1^2 - \frac{3}{4}\gamma z_1^3 + F(d) = \sigma(z_1) = 0, \end{cases} \quad (13)$$

where $\sigma(z_1)$ represents the restoring force acting on the system and plays a vital role in the existing bifurcations in the dynamics of the system. Calculating the first derivative of $\sigma(z_1)$ with respect to $z_1$ produces the corresponding stiffness function

$$K(z_1) = \frac{d\sigma}{dz_1} = -J(d) + 2G(d)z_1 - \frac{9}{4}\gamma z_1^2, \quad (14)$$

for the system, which reveals how the stiffness function changes with $d$ and $z_1$. The corresponding potential function is obtained by taking integral of terms associated with $z_1$ according to Eq. (20)

$$V(x) = J(d)\frac{z_1^2}{2} - G(d)\frac{z_1^3}{3} + \frac{3}{16}\gamma z_1^4, \qquad (20)$$

The potential function is plotted in Figure 3a for $\gamma = 150$ as well as its shape for three values of slow variable $d$ is depicted in the Figure 3b. These figures reveal that the potential function changes from double well to single well as slow parameter varies, such that the system is continuously experiencing qualitative change as time passes. Consequently, the equilibrium points of the system not only change significantly, they also change in an oscillatory manner due to the time dependance of the slow parameter $d$.

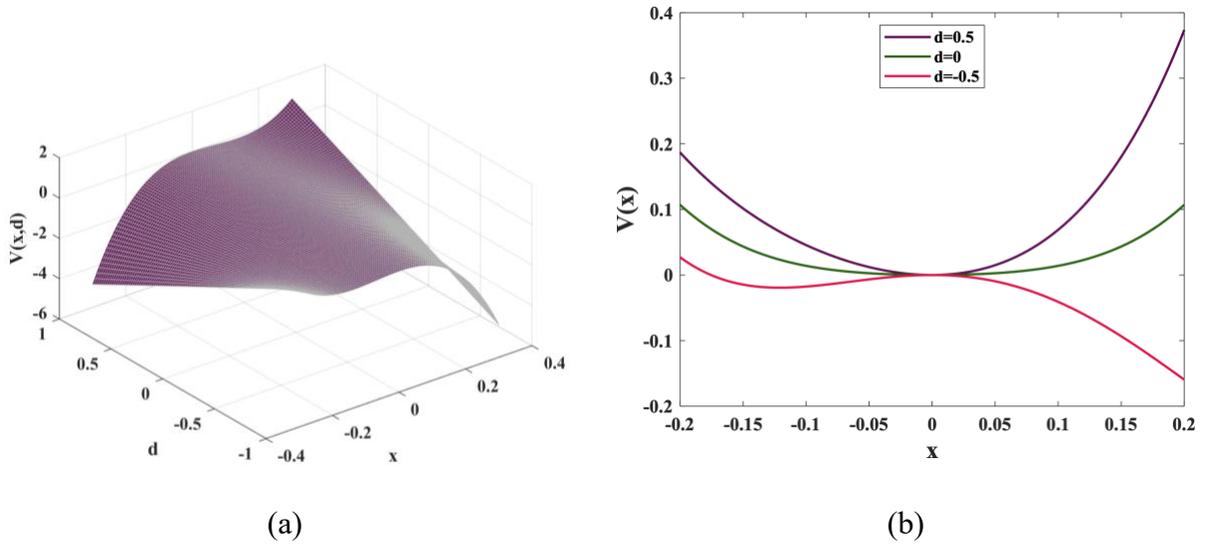

(a)          (b)

Figure 3. For $\gamma = 150$ (a) Potential function variation with $d$ and $x$. (b) The potential function variation from single well to double well by changing slow variable $d$.

To perform the local bifurcation analysis, we first compute the Jacobian matrix of the autonomous system Eq. (14) at the equilibrium point $z_1 = z_1^*$, and obtain

$$A = \begin{pmatrix} 0 & 1 \\ -K(z_1^*) & -\xi \end{pmatrix}, \qquad (21)$$

where $\xi$ is damping ratio. The characteristic equation at the equilibrium point is then

$$\lambda^2 + \xi\lambda + K(z_1^*) = 0, \qquad (22)$$

which has eigenvalues $\lambda_{1,2} = \mu \pm i\psi$, where

$$\mu = \frac{-\xi}{2}, \quad \psi = \sqrt{\xi^2 - 4K(x)}/2. \tag{15}$$

According to Ruth-Hurwitz stability, the coefficients of the Eq. (22) must be positive for the system to be stable. For $\xi^2 = 4K(x)$, the system's eigenvalues are $\lambda_{1,2} = -\xi/2$ and response is stable. For $\xi^2 < 4K(x)$, bifurcations such as Hopf does not occur because the real component is nonzero. For $4K(x) = 0$, then $\lambda_1 = 0$ and $\lambda_2 \neq 0$, therefore fold bifurcation occurs in the dynamics, cusp bifurcation get the possibility to occur which is a higher-order degeneracy of fold (or limit point) bifurcation [30]. Here, we discuss fold and its mathematical definition.

- **Definition of fold bifurcation [30]:** suppose the system

$$\dot{x} = f(x, \alpha), \quad x \in \mathbb{R}^1, \quad \alpha \in \mathbb{R}^1, \tag{16}$$

with a smooth $f$ has the equilibrium $x = 0$ with $\lambda = f_x(0,0) = 0$ at $\alpha = 0$, where $f_x$ is first derivative of $f$ with respect to $x$. Expanding $f(x, \alpha)$ using a Taylor series with respect to $x$ at $x = 0$:

$$f(x, \alpha) = f_0(\alpha) + f_1(\alpha)x + f_2(\alpha)x^2 + o(x^3). \tag{17}$$

Two conditions are satisfied: $f_0(0) = f(0,0) = 0$ (equilibrium condition) and $f_1(0) = f_x(0,0) = 0$ (fold bifurcation condition). The $f_0, f_1$ and $f_2$ are coefficients obtained by Taylor series expansion.

In conjunction with the provided definition of fold, it is useful to consider that the normal form of fold bifurcation of a scalar system is

$$\dot{z}_1 = \alpha + z_1^2 + o(z_1^3), \tag{18}$$

which is locally topologically equivalent near the origin to the system $\dot{z}_1 = \alpha + z_1^2$ [30]. This is because adding the higher order terms does not affect the fold bifurcation at neighborhood of equilibrium point [30]. Making analogy between Eq. (14) and Eq. (26) explains coefficient of the term $q(t)$ play role in folding points.

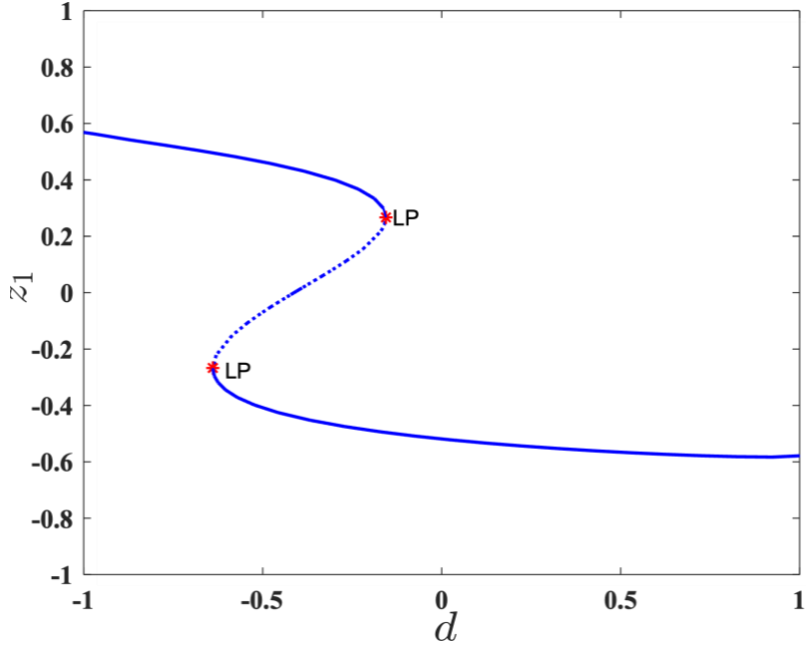

Figure 4. One parameter bifurcation diagram of fast subsystem in terms of slow variable $d$ for $\eta_0 = 0.1$, $\sigma = 120, \gamma = 150$. Red markers (LP) show fold bifurcations, (....) unstable equilibrium.

Hereunder, we investigate the system behavior at the fold bifurcation numerically. The usefulness of performing bifurcation analysis lies in partitioning the parameter space of a problem into distinct regions, where the dynamical system is expected to exhibit qualitatively distinct behavior. Bifurcation diagrams in this work are generated using the MATCONT toolbox, which is based on continuation methods. The process begins by finding the equilibrium points of the system, followed by analyzing the system's stability for each parameter value through eigenvalue calculations. Two generic codimension-1 bifurcations that can be detected are the fold and the Hopf bifurcations. Every value of a parameter that leads to one zero eigenvalues, stored as a fold point. With the aim of obtaining a good understanding how elastic foundation and excitation of a traveling wave affects the dynamic response of a beam, the foundation properties $(\sigma, \gamma)$ as well as traveling wave parameters $(\kappa, \eta_0, \omega)$ are considered as bifurcation parameters throughout this work. However, the focus of this study is on fold bifurcations as opening window for cusp and Bogdanov-Takens bifurcations. To investigate the instabilities of the system, the parameters are set at $\eta_0 = 0.1$, $\xi = 0.02$ which are appropriate values for the investigated beam-foundation configuration. The manifold of equilibria or critical manifold is shown in the Figure *4* which is obtained by setting fast subsystem $\dot{z}_1 = 0$ in Eq. (14). Fold points are also shown with red markers (LP or limit point) as well as the unstable segments of critical manifold is plotted with dotted line. Fold bifurcations arise when stable and unstable equilibria interact and annihilate each other. Hence, existence of pairs of fold bifurcations often lead to switching or jumping behavior at vicinity of folding points, which is where the response jumps from one equilibrium to another. This happens when a region of bi-stability exist. In fact, for the system studied here, the external excitation forces it to jump from one equilibrium to another instantaneously which is shown in the Figure 5.

Figure 5, for two values of $\gamma$ which represents contribution of nonlinear foundation illustrates the phase portrait and corresponding transformed version of it of the mixed-mode oscillation is presented and the curve of equilibrium points, shown with dashed line is also overlayed and reveals that the how response pattern of the system drastically change by approaching to the unstable manifold. The numerical simulation shows reasonable agreement between the transformed phase diagram and the bifurcation diagram which complements the fact that it is the switching between the two stable/unstable equilibria that form the mixed mode oscillation in this system. When the system reaches folding points on manifold of equilibria in which one has attracting and one has repelling character, it loses its stability and brings the onset of bursting responding. In other words, the limiting behavior of the system is jumping to another state which is known by a fast transition. Thus, the system switching back and forth between stable and unstable points. In fact, adding a nonlinear foundation in conjunction with base excitation to the linear beam leads to a bi-well potential function. Also, the isolated regions in transformed phase portraits correspond to stable area therefore no jumping behavior should be expected.

A key question arises here: for what combination of linear and nonlinear components of the foundation stiffness leads to a fold bifurcation, which marks the onset of instability in the system? This is addressed through the 2D bifurcation diagram shown in the Figure 7a which exhibits for each value of linear stiffness ($\sigma$) at two values of nonlinear stiffness ($\gamma$) fold occurs. As could be seen, stiffer linear foundation demands higher nonlinear contribution to imposing the folding behavior. Increasing trend in the both parameters signifies that folding occurs at stiffer foundation which by comparing the overall trend depicted in the Figure 7b, it could be concluded that by increasing stiffness, folding points move toward to the neutral position of the beam ($z_1 =$ 0), which means at shorter distances from the neutral position the beam is able to fluctuates, thus the beam is less prone to instability.

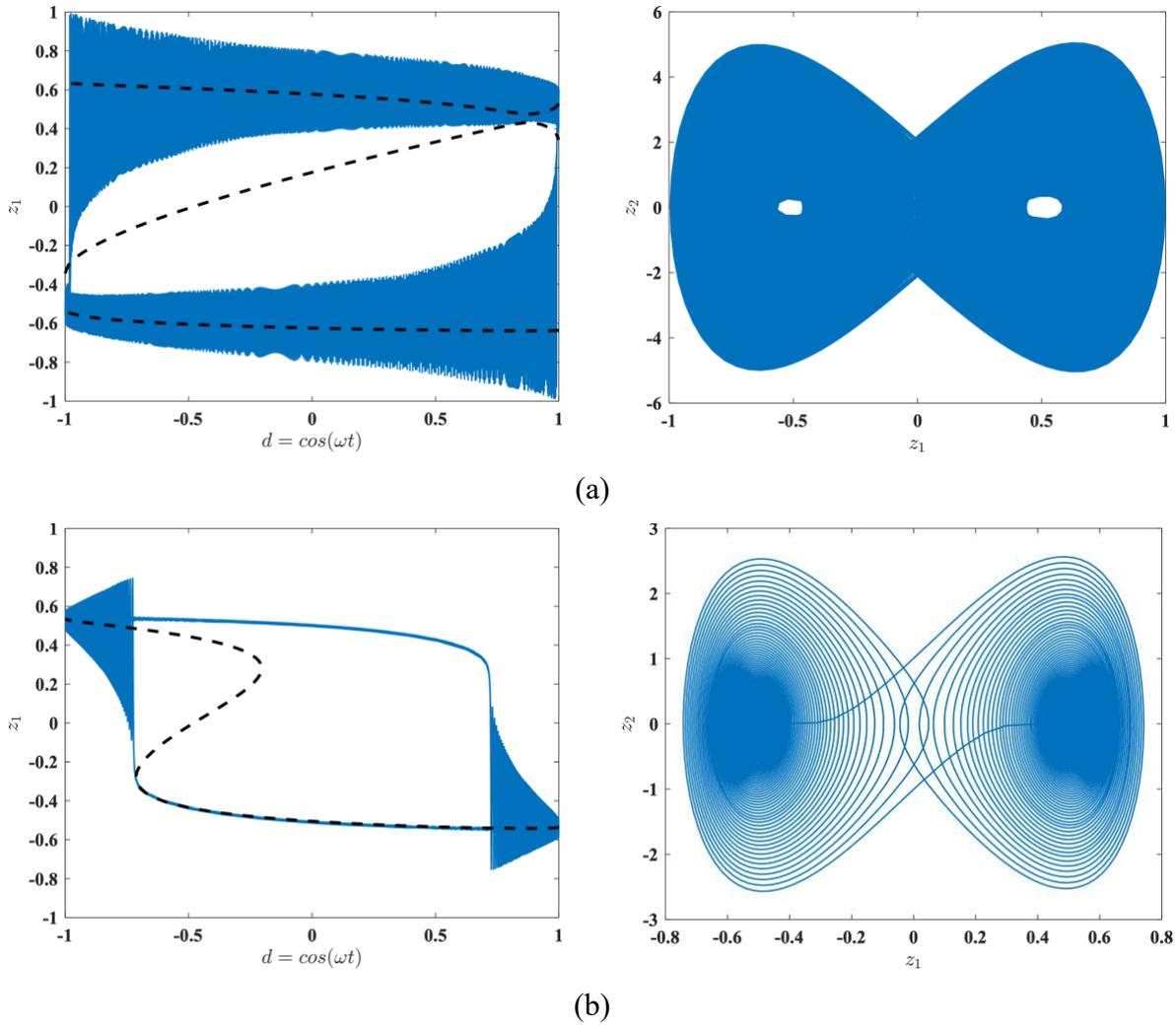

Figure 5. Bursting oscillations for transformed phase portrait and phase portrait overlaid with the equilibrium manifold (dashed black line) for $\omega = 0.01, \eta_0 = 0.1$, (a) $\gamma = 12$, (b) $\gamma = 150$.

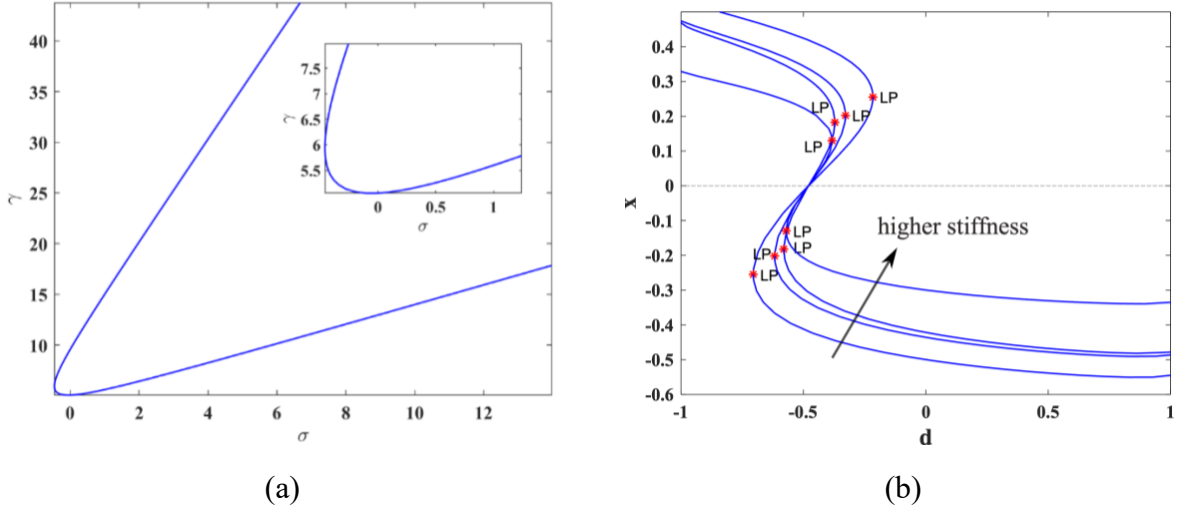

(a)                                     (b)

Figure 6. (a) 2D bifurcation map of folding point in parameters space $(\sigma, \gamma)$, (b) equilibrium manifold vs stiffness of foundation.

Considering the importance of amplitude of excitation in the bursting behaviors, the bifurcation structure of system Eq. (9) on the space $(d, \eta_0)$ using parameter sweeping and continuation method for constant $\sigma$ and $\gamma$ is depicted in Figure 7a which separate the whole parameter space by two curves. The lower curves show that by evolving slow variable to near zero and positive values, the fold occurs at very lower amplitudes of the wave. This result is interrelated with trend showed in Figure 6b, since higher wave amplitudes leads to higher stiffness of the system and higher stiffness is less prone to fold behavior. In summary, wave amplitude $\eta_0 < 0.1$ is more critical working region. The upper curve manifests very high amplitudes which has no physical importance and they are not consistent with the assumptions of the Euler-Bernoulli beam theory used in this work. Parallel to wave amplitude, for investigating effect of wave number effect on the system qualitative changes, 2D bifurcation study is performed on parameter space $d, \kappa)$, see Figure 8b. According to the region enclosed by the blue lines, for major range of $\kappa$ system passes through a fold point two times by evolving slow variable, the result that is expected from two-fold points in the Figure 4. By going to higher wave numbers, the curves become narrower which means shorter space is available to response under effect of fold bifurcation. It is noteworthy to mention that the wave number which controls spatial frequency of the wave is proportional in inverse of wavelength $(H)$, $\kappa = 2\pi/H$. Thus, low $\kappa$ value corresponds to base excitation which changes slowly over the length of the beam and high $\kappa$ value means the beam is excited by a base subjected to rapid curvy wave profile with multiple peaks over the length of the beam. The red markers imply

Bogdanov-Takens (BT) bifurcations, also known as double zero eigenvalues, since at BT points determinant and Jacobian of system of Eq. (14) becomes zero. The BT points may act as intersection or origination center for emerging other bifurcations such as Hopf or homoclinic. The goal here is just to show that the system has another kind of bifurcations than fold and cusp as well as showing that all parameters of the elastic foundation and traveling wave has potential to create bifurcation, which could be considered for future research.

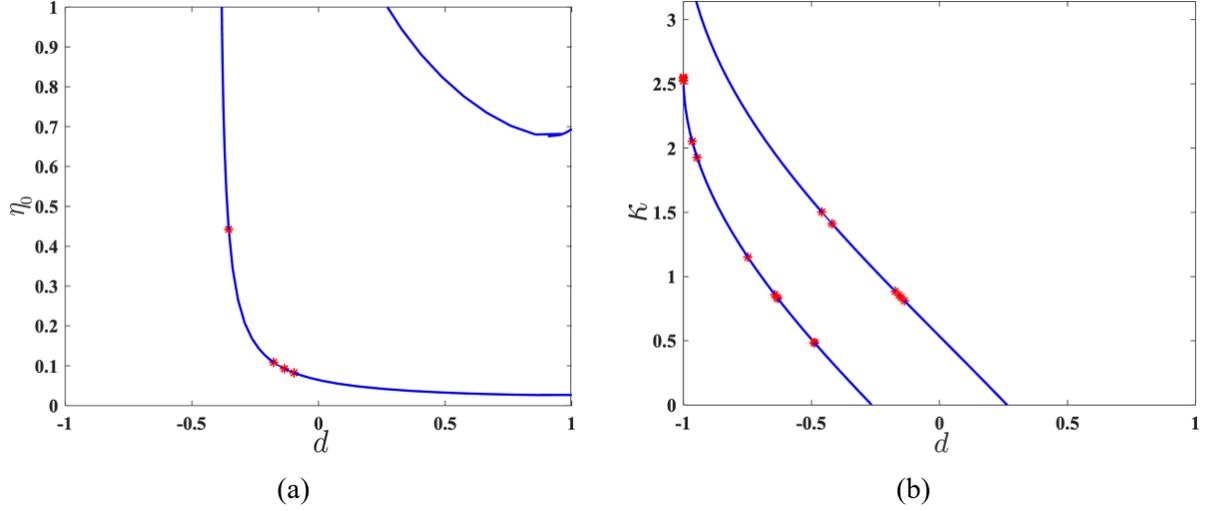

(a) (b)

Figure 7. 2D bifurcation set of the fast subsystem. (—) fold bifurcation, (*) Bogdanov-Takens bifurcations.

Usually, in nonlinear oscillations mainly based on perturbation methods, capability of methods are affected by type of nonlinearity (strong or weak), but in slow-fast dissection method there is not such a limitation therefore it allows to study nonlinear systems in unified manner via exploring the entire phase space of the system. However, this method has sensitivity to perturbation in parameters as well as admitting only integer multipliers of frequencies.

The definition of cusp bifurcation and conditions of occurrence are provided here.

- **Definition of cusp bifurcation [31]:** Let $f \in C^3(\mathbb{R} \times \mathbb{R}^3, \mathbb{R})$ and

$$f(0,0) = D_x f(0,0) = 0, \qquad (19)$$

$$D_x^2 f(0;0) = 0, \qquad (20)$$

$$D_x^3 f(0;0) \neq 0, \qquad (21)$$

Let $N$ be the neighborhood of $\mu = 0$ for which $x_0(\mu)$ is the unique solution of $D_x^2 f(x_0(\mu); \mu) = 0$ such that $x_0(0) = 0$. Then $f$ has a cusp bifurcation in $N$ with bifurcation set $27A^2(\mu)D(\mu) = -4B^3(\mu)$ where

$$A(\mu) = f(x_0; \mu), \quad B(\mu) = D_x f(x_0; \mu), \quad D(\mu) = \tfrac{1}{6} D_{xxx} f(x_0; \mu).$$

Results discussed so far used stiffness of foundation at the order of the beam stiffness ($EI$). Considering time variable contribution of stiffness, with the purpose of exploring another instabilities the lower stiffness parameters were used as $\sigma = 12$ and $\gamma = 1.5$ which has led to occurrence of cusp bifurcation as illustrated in the Figure 8. Near a cusp point, number of steady states on manifold of equilibria may vary from one to three depends on $(d, \kappa)$ parameters. In case of three steady state coexistence, the middle point is unstable therefore, abrupt transition at two-fold bifurcations to the left or to the right of manifold of equilibria may be expected.

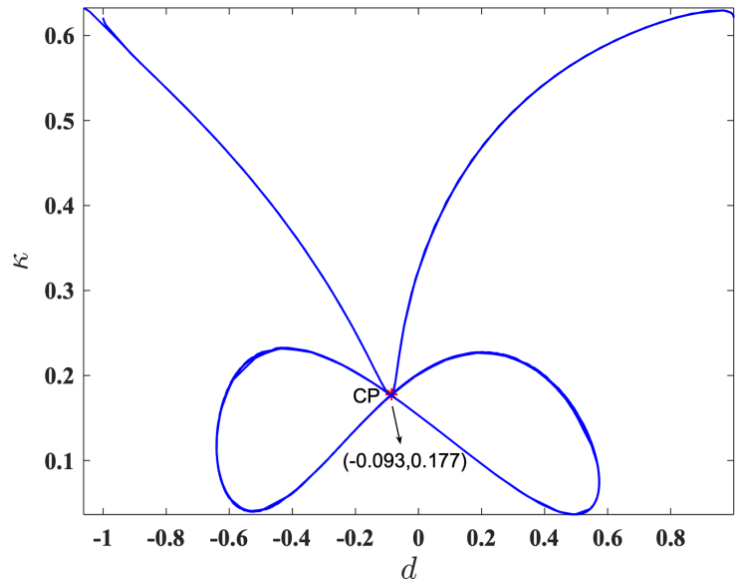

Figure 8. Cusp bifurcations caused on $(d, \kappa)$ parameter space

## 4. Concluding Remarks

Bursting behavior is one of the complex types and investigating them is vital in multiscale dynamical systems. In this paper we present bursting patterns caused by low frequency excitations for a nonlinear beam-elastic foundation under simultaneous parametric and external multi-frequency excitation caused by traveling wave using fast-slow technique. Main idea behind this work is to show orders gap difference between natural frequency of the system and excitation frequency leads to less reported oscillation patterns including bursting pattern. Two-dimensional bifurcation maps in conjunction with transformed phase diagrams showed that linear stiffness, amplitude of excitation and frequency plays a key role in mixed mode oscillation. It is shown that the nonlinear beam-elastic foundation configurations are capable to exhibit diverse oscillatory responses and multiplicity of excitation frequencies adds to the complexity of system response through adding more fold points to the equilibria manifold. In this work it was assumed that there are exact frequencies for excitation, while in practice they can perturb. So, more work could be done to explore conditions of perturbation of excitation frequencies. Also, considering the applications of beam and elastic foundation in modeling, the novel findings in this works are not only insightful for the cubic nonlinearity discussed in this work, but also similar instabilities may be found in other types of structure and elastic foundations interactions. This work shows exerting a base excitation to a beam using a traveling wave, adding complex parametric multi-frequency excitations and alter the dynamics of the system substantially; a situation which does not happen for beam resting on elastic foundation under concentrated or distributed loadings. On the other hand, although the nonlinear foundation imposes instabilities to the system mainly in the form of folding bifurcations, it brings a self-stabilizing effect to the system.

## Funding Declaration

The authors declare that no external funding is used for this research. The authors are grateful for internal start-up funding provided by the Georgia Institute of Technology that supported this work.

## Competing interests

The authors have no relevant financial or non-financial interests to disclose.

## Author Contributions

**S.M:** Writing-Original Draft, Methodology, Software, Validation, Formal analysis, Investigation, Data curation, Conceptualization. **K.J.M:** Writing-Original Draft, Methodology, Software, Validation, Formal analysis, Supervision, Investigation, Data curation, Conceptualization, Project administration.